\documentclass[12pt]{article}%
\usepackage{graphicx}
\usepackage[intlimits]{amsmath}
\usepackage{latexsym}
\usepackage{amsfonts}
\usepackage{amssymb}%
\makeatletter
\newfont{\footsc}{cmcsc10 at 8truept}
\newfont{\footbf}{cmbx10 at 8truept}
\newfont{\footrm}{cmr10 at 10truept}
\newtheorem{theorem}{Theorem}

\begin{document}

\title{The Moore bound for Spectral Radius}
\author{Vladimir Nikiforov\\{\small Department of Mathematical Sciences, University of Memphis, }\\{\small Memphis TN 38152, USA}}
\maketitle

\begin{abstract}
Let $G$ be a graph with $n$ vertices, $m$ edges, girth $g,$ and spectral
radius $\mu.$ Then
\[
\mu^{2}\leq\left\{
\begin{array}
[c]{ll}%
m & \text{if }g=4\\
n-1 & \text{if }g\geq5.
\end{array}
\right.
\]
\medskip

\textbf{Keywords: }\medskip Moore bound, spectral radius, irregular graph,
degree sequence

\end{abstract}

What is the maximum spectral radius of a graph of order $n$ and girth $g$?
Equivalently, what is the minimum order of a graph with girth $g$ and spectral
radius $\mu.$ For $d$-regular graphs the spectral radius is equal to $d;$
thus, the Moore bound (see, e.g., \cite{Big93}, p.180 and \cite{MiSi05} for a
general survey) shows that the order of such graphs is at least
\[
n\left(  d,g\right)  =\left\{
\begin{array}
[c]{ll}%
1+d+d\left(  d-1\right)  +\cdots+d\left(  d-1\right)  ^{r-1} & \text{if
}g=2r+1\\
2\left(  1+\left(  d-1\right)  +\cdots+\left(  d-1\right)  ^{r-1}\right)  &
\text{if }g=2r.
\end{array}
\right.
\]

Alon, Hoory, and Linial \cite{AHL02} extended this result, showing that every
graph with girth $g$ and \emph{average degree} $d$ has at least $n\left(
d,g\right)  $ vertices.

The bound changes radically when the average degree is replaced by the
spectral radius.

\begin{theorem}
\label{th1}Let $G$ be a graph with $n$ vertices, $m$ edges, girth $g,$ and
spectral radius $\mu.$ Then%
\[
\mu^{2}\leq\left\{
\begin{array}
[c]{ll}%
m & \text{if }g=4\\
n-1 & \text{if }g\geq5.
\end{array}
\right.
\]
If $g=4,$ then $\mu^{2}=m$ holds if and only if $G$ is a complete bipartite
graph possibly with isolated vertices. If $g\geq5$ then $\mu^{2}=n-1$ holds if
and only if $G$ is either a Moore graph of diameter $2$ or a star.
\end{theorem}

For $g=4$ this result was obtained in \cite{Nos70} and extended in another
direction in \cite{Nik02}. The case $g\geq5$ was given in \cite{FMS93}.


\begin{thebibliography}{9}                                                                                                %


\bibitem {AHL02}N. Alon, S. Hoory, N. Linial, The Moore bound for irregular
graphs, \emph{Graphs Combin.} \textbf{18} (2002), 53--57.

\bibitem {Big93}N. Biggs, \emph{Algebraic graph theory,} Cambridge University
Press, Cambridge, second edition, 1993.

\bibitem {FMS93}O. Favaron, M. Mah\'{e}o, J.-F. Sacl\'{e}, Some eigenvalue
properties in graphs (conjectures of Graffiti. II), \emph{Discrete Math.}
\textbf{111} (1993), 197--220.

\bibitem {MiSi05}M. Miller, J. \v{S}ir\'{a}\v{n}, Moore graphs and beyond: A
survey of the degree/diameter problem,
\emph{www.combinatorics.org/Surveys/ds14.pdf.}

\bibitem {Nik02}V. Nikiforov, Some inequalities for the largest eigenvalue of
a graph, \emph{Combin. Probab. Comput}. \textbf{11} (2002), 179--189.

\bibitem {Nos70}E. Nosal, Eigenvalues of Graphs, Master's thesis, University
of Calgary, 1970.
\end{thebibliography}
\end{document}